\documentclass[preprint,3p,times]{elsarticle}
\usepackage{prettyref}
\usepackage{url}
\usepackage{amsthm}
\usepackage{amsmath}
\usepackage{graphicx}
\usepackage{amssymb}
\usepackage{esint}

  \theoremstyle{remark}
  \newtheorem*{rem*}{Remark}
\theoremstyle{plain}
\newtheorem{thm}{Theorem}
  \theoremstyle{plain}
  \newtheorem{lem}[thm]{Lemma}

\makeatletter
\def\ps@pprintTitle{%
  \let\@oddhead\@empty
  \let\@evenhead\@empty
  \def\@oddfoot{\reset@font\hfil\thepage\hfil}
  \let\@evenfoot\@oddfoot
}
\makeatother
\biboptions{sort&compress}

\begin{document}
\begin{frontmatter}

\title{On the noise modelling in a nerve fiber}

\author[UoL,fti]{A.A. Samoletov\corref{cor1}}

\ead{samolet@fti.dn.ua}
\ead[url]{http://fti.dn.ua/~samolet}

\author[UoL]{B.N. Vasiev\corref{cor2}}

\ead{b.vasiev@liverpool.ac.uk}
\ead[url]{http://www.maths.liv.ac.uk/~bnvasiev}

\cortext[cor2]{Corresponding author}

\cortext[cor1]{Principal corresponding author}

\address[UoL]{Department of Mathematical Sciences, University of Liverpool, Peach
Street, Liverpool, L69 7ZL, UK}

\address[fti]{Institute for Physics and Technology, NASU, 72 Luxsembourg Street,
83114 Donetsk, Ukraine}

\begin{abstract}
We present a novel mathematical approach to model noise in dynamical
systems. We do so by considering dynamics of a chain of diffusively
coupled Nagumo cells affected by noise. We show that the noise in
transmembrane current can be effectively modelled as fluctuations
in electric characteristics of the membrane. The proposed approach
to model noise in a nerve fibre is different from the standard additive
stochastic current perturbation (the Langevin type equations).\end{abstract}

\end{frontmatter}

\section{Introduction}

A typical nerve fiber is coated in myelin (the myelin sheath consists
of a single Schwann cell which is wrapped about $100$ around the
nerve fiber) with spatially periodic gaps, the nodes of Ranvier. Roughly
the myelin sheath increases the membrane resistance by a factor of
about $100$ and decreases the membrane capacitance by a factor of
about $100$. Typically, the width of the node of Ranvier is about
$1\mathrm{\mu m}$, the distance between nodes (the length of myelin
sheath) is about $1.5\mathrm{mm}$ that is close to $100d$, where
$d$ is the nerve fiber diameter. Transmembrane ion flow occurs only
at the nodes of Ranvier. Diffusive coupling corresponds to axial currents
between nodes and allows propagation of changes in the transmembrane
potential (action potential) in the spatial variable. Nerve fibers
behave as intrinsic spatially discrete systems. The biological reason
for such discrete structure: propagation of action potential along
myelinated fibre is faster compare to that in nonmyelinated because
of its saltatory propagation between nodes (speed in a myelinated
fibre is around$100\, m/s$ and in a nonmyelinated fiber is $1\div5\, m/s$).
For further details of the model, physical parameters and equivalent
electric circuit we refer to \citep{Keener2009}. 

The idea that noise can play a positive role and benefit neural function
is relatively new. Just 60 years ago it was commonly accepted that
the noise is destructive to neural encoding \citep{Fatt1950,Fatt1952}.
Today it is well established that noise plays constructive role in
the nerve system \citep{Zucker2003,Sharma2003,Milton2005,Faisal2008,Deco2009}.
This new paradigm was initiated by research on stochastic resonance
phenomena. It was shown that the stochastic resonance improves the
transfer of information \citep{Gammaitoni1998,Samoletov2004}. 

In this Letter, we accomplish two goals. First, we propose and study
the deterministic scheme for modelling of noise in a nerve fiber. This
scheme involves dynamical fluctuations of electric characteristics
of the membrane together with their negative feedback control depending
on the noise intensity. Then, to ensure ergodicity property of the
dynamics, we combine this dynamical feedback control with a stochastic
perturbation. In contrast to the random noise model (see \prettyref{sec:Random-noise})
our scheme operates with the only white noise process that indirectly
affect initial dynamics. While we do not claim that our scheme to
model noise in a nerve fibre is better than standard additive stochastic
current perturbation (the Langevin type equations), we state our approach
as different.

We consider a lattice of diffusively coupled Nagumo cells described,
in absence of noise, by the equations, $\begin{array}{l}
\dot{u_{i}}=l\triangle u_{i}+f(u_{i}),\:\end{array}i\in\mathbb{Z}$ is a spacial index, where $f(u):\mathbb{R}\rightarrow\mathbb{R}$
has a bistable character,\emph{ }for example $f(u)=-ku(u-\alpha)(u-1)$,
$0<\alpha<1$, $k>0$; $\Delta u_{i}\equiv u_{i+1}-2u_{i}+u_{i-1}$
is the standard $3$ point discretization of the Laplacian (discrete
Laplacian), and $l>0$ is a coefficient of the diffusive coupling.
In addition we define {}``potential'' $V(u)$ by the differential
equation, $V'(u)=-f(u),\: V(0)=0$. In these equations variable $u$
corresponds to a transmembrane electric potential, $k$ corresponds
to the membrane conductance, $\alpha$ is the threshold potential.
Besides modelling of the action potential propagation along a nerve
fibre, this lattice system is important in another different areas
of research \citep{Cahn1960,Erneux1993,Bates1999}.

In cases where $i\in I\subset\mathbb{Z}$ and $I$ is bounded, we
consider this set with respect to boundary conditions, for example
of the Neumann type. For what follows, it is convenient to represent
Nagumo equations in the variational form. Define the {}``energy''
functional, $\mathcal{V}[u]=\sum_{\{i\}}[\frac{1}{2}l(\nabla u_{i})^{2}+V(u_{i})]$,
where $\nabla u_{i}=u_{i}-u_{i-1}$ is the discrete gradient. Hereafter
we accept short notations: $\frac{\partial}{\partial u_{i}}\equiv\partial_{i}$,
$\frac{\partial^{2}}{\partial u_{i}^{2}}\equiv\partial_{i}^{2}$,
$\frac{\partial}{\partial t}\equiv\partial_{t}$, and so on. With
these definitions, rewrite the lattice of the diffusively coupled
Nagumo equations in the gradient form,\begin{equation}
\dot{u}_{i}=-\partial_{i}\mathcal{V}[u],\quad i\in\mathbb{Z}.\label{eq:N}\end{equation}
It is easy to reveal that $\mathcal{V}[u]$ is the Lyapunov functional
since $\dot{\mathcal{V}}[u]=-\sum_{\{i\}}\left(\partial_{i}\mathcal{V}[u]\right)^{2}\leq0$.
Steady states of equations \prettyref{eq:N} are the extrema of functional
$\mathcal{V}[u]$. The minima and maxima of $\mathcal{V}[u]$ correspond
respectively to stable and unstable solutions of equation \prettyref{eq:N}.
Suppose $\mathcal{V}[u]>-\infty$ and define for a continuous function
$A(u)$ the time averaging, $\overline{A(u)}=\lim_{T\rightarrow\infty}\frac{1}{T}\int_{0}^{T}A(u(t))dt$.
Applying the time averaging to $\dot{\mathcal{V}}[u]=-\sum_{\{i\}}\left(\partial_{i}\mathcal{V}[u]\right)^{2}$
($\mathcal{V}[u]>-\infty$) we arrive at the equation, $\overline{\sum_{\{i\}}\left(\partial_{i}\mathcal{V}[u]\right)^{2}}=0.$
Thus the system spends almost all time at states of extrema of $\mathcal{V}[u]$.
These extrema are solutions of the discrete lattice equation, $l\triangle u_{i}-\partial_{i}\mathcal{V}[u]=0,\quad i\in\mathbb{Z}$.
This equation is implied to be equipped with some boundary conditions.
For our purpose we accept the following conditions, $\nabla u_{i}\rightarrow0\quad\mathrm{as}\quad i\rightarrow\pm\infty.$

\section{Random noise\label{sec:Random-noise}}

To model the influence of noise on deterministic system (\ref{eq:N}),
it is widely accepted in the literature that the noise is implemented
in (\ref{eq:N}) by the additive stochastic currents, $\xi_{i}(t),\, i\in\mathbb{Z}$,
where $\left\{ \xi_{i}(t)\right\} _{i\in\mathbb{Z}}$ is the set of
independent standard generalized Gaussian $\delta$-correlated processes
completely characterized by the first two cumulants, $\left\langle \xi_{i}(t)\right\rangle =0$
and $\left\langle \xi_{i}(t)\xi_{j}(t^{\prime})\right\rangle =\delta_{ij}\delta(t-t^{\prime})$;
$<...>$ means averaging over all realizations of the random perturbations.
The set of stochastic differential equations corresponding to \prettyref{eq:N}
takes the form, $\begin{array}{l}
\dot{u_{i}}=l\triangle u_{i}+f(u_{i})+\sqrt{2D}\xi_{i}(t)\end{array}=-\partial_{i}\mathcal{V}[u]+\sqrt{2D}\xi_{i}(t),\quad i\in\mathbb{Z}$, where $D$ is the noise intensity (we suppose that noise does not
depend on node). It is convenient for what follows to represent this
set as \begin{equation}
\begin{array}{l}
\dot{u_{i}}=\end{array}-\lambda\partial_{i}\mathcal{V}[u]+\sqrt{2\lambda D}\xi_{i}(t),\quad i\in\mathbb{Z},\label{eq:NS}\end{equation}
where a reference time scale $\lambda$ is explicitly introduced.
Rescaling time in \prettyref{eq:NS}, $t\rightarrow\lambda^{-1}t$,
and taking into account scaling property of the white noise, we arrive
at the case $\lambda=1$.

In system \prettyref{eq:NS} dissipative processes and random perturbations
equilibrate one another. In respect of the {}``energy'', $\mathfrak{\mathcal{V}}[u]$,
we arrive at the stochastic differential equation (we specify this
equation in the sense of Stratonovich (\emph{e.g.} \citep{Oksendal2003})),
\begin{equation}
\mathcal{\dot{V}}[u]=-\lambda\sum_{\{i\}}\left(\partial_{i}\mathcal{V}[u]\right)^{2}+\sqrt{2D\lambda}\sum_{\{i\}}\partial_{i}\mathcal{V}[u]\xi_{i}(t).\label{eq:SLyap}\end{equation}
Equation \prettyref{eq:SLyap} demonstrates in what way the noise
affects the {}``energy'' and defines rate of its (stochastic) fluctuations.
Assume that $\mathcal{V}[u]>-\infty$. Then after averaging over all
realization of the random perturbations we arrive at the relation,
$-\left\langle \sum_{\{i\}}\left(\partial_{i}\mathcal{V}[u]\right)^{2}\right\rangle +D\left\langle \sum_{\{i\}}\partial_{i}^{2}\mathcal{V}[u]\right\rangle =0$,
that does not depend on $\lambda$; we assume $\left\langle V[u]\right\rangle =const$.
The relation can be derived either by an elementary calculation or
elegantly applying Novikov's formula \citep{Klyatskin2005a}. This
is an important relation that connects the noise intensity to configurational
ensemble averages and thus can be considered as the definition of
the noise intensity. In what follows we conjecture that the analogue
formula involving the time averaging instead of the ensemble averaging,
$-\overline{\sum_{\{i\}}\left(\partial_{i}\mathcal{V}[u]\right)^{2}}+D\overline{\sum_{\{i\}}\partial_{i}^{2}\mathcal{V}[u]}=0$,
is valid and thus defines the noise intensity in the framework of
deterministic dynamics. In order to deepen the conjecture and to describe
dynamics of the deterministic fluctuations, we have to further presume
the rate of dynamic fluctuations (r.d.f.) in the form, \begin{equation}
\mathrm{r.d.f.}\sim-\sum_{\{i\}}\left(\partial_{i}\mathcal{V}[u]\right)^{2}+D\sum_{\{i\}}\partial_{i}^{2}\mathcal{V}[u],\label{eq:Dyn-D-rate}\end{equation}
that is instead of random perturbations (that do not present in deterministic
dynamics) we need to consider dynamic fluctuations of an appropriate
variable. Indeed, in absence of random perturbations we have to adopt
an another way to properly perturb the system. Fluctuations in the
electric characteristics of membrane are conjugate to that in electric
current across the membrane. Thus it is reasonable to consider a certain
electric characteristic of the membrane, supposedly $RC$ (where $R$
is resistance and $C$ is capacitance), that defines a time scale
and allows this characteristic to dynamically fluctuate. 

With random noise its intensity $D$ is commonly considered as an
independent parameter. Indeed, the Fokker-Planck operator corresponding
to \prettyref{eq:NS} has the form, $\mathcal{F}^{*}\rho\equiv-\sum_{(i)}\partial_{i}\left(\partial_{i}\mathcal{V}[u]\rho\right)+D\sum_{(i)}\partial_{i}^{2}\rho$.
The Fokker-Planck equation associated with $\mathcal{F}^{*}$, $\partial_{t}\rho=\mathcal{F}^{*}\rho$,
allows the invariant solution, $\rho_{\infty}\left[u\right]\sim\exp\left\{ -D^{-1}\mathcal{V}[u]\right\} $.
We prove the identity, $\mathcal{F}^{*}\rho_{\infty}\left[u\right]\equiv0$,
by straightforward calculation. It is known that this distribution
and the corresponding probabilistic measure, $d\mu\sim\exp\left\{ -D^{-1}\mathcal{V}[u]\right\} \prod_{(i)}du_{i}$,
are typically unique for dynamics \prettyref{eq:NS}. In a word, stochastic
dynamics \prettyref{eq:NS} is typically ergodic. This means that
for every continuous function $A$, $\int A(u)d\mu=\lim_{T\rightarrow\infty}\frac{1}{T}\int_{0}^{T}A(u(t))dt$,
almost for sure for all initial values $u(0)$. The invariant measure
relates to infinite time interval. Thus scaling the time variable
does not affect the measure. The invariant (equilibrium) distribution
$\rho_{\infty}$ demonstrates explicit dependence on the noise intensity,
$D$. The only constraint on $D$ arise when we presume a nondestructive
role of the noise. Namely, in the case of a cubic nonlinearity of
$f(u)$, the general form of $V(u)$ is double-well. Then the noise
can induce transition from one well to another, it depends on $D$,
and is expected to be a slow process.

Now we can pose the problem: Given a probability measure $d\mu$ (or
an augmented measure on an extended phase space). It is necessary
to find a dynamics such that $\lim_{T\rightarrow\infty}\frac{1}{T}\int_{0}^{T}A(u(t))dt=\int A(u)d\mu$
for every continuous function $A$. We say this is dynamic modelling
of the noise and assume $d\mu\sim\exp\left\{ -D^{-1}\mathcal{V}[u]\right\} \prod_{(i)}du_{i}$
as the invariant (ergodic) measure for this dynamics.

\section{Deterministic modelling of noise}

Now we will put together the above observations, - that the rate of
feedback control of dynamic fluctuations and the invariant measure,
depend on the noise intensity, - to derive a model of a deterministic
noise of intensity $D$ in a nerve fibre. The requirements are:
\begin{itemize}
\item Dynamics of $u$ depends on the external dynamic variables (\emph{e.g.},
$\lambda$ is endowed with its own equation of motion);
\item Rate of deterministic dynamic fluctuations is directly related to
\prettyref{eq:Dyn-D-rate} (\emph{e.g.}, the rate of fluctuations
is a measure of the influence of environment on electrical characteristics
of the membrane);
\item Measure $d\mu\sim\exp\left\{ -D^{-1}\mathcal{V}[u]\right\} \prod_{(i)}du_{i}$
is invariant for the dynamics;
\item Dynamics is ergodic.
\end{itemize}
In other words, we will sample the invariant measure, $d\mu\sim\exp\left\{ -D^{-1}\mathcal{V}[u]\right\} \prod_{(i)}du_{i}$,
by the method proposed in \citep{SDC07,Samoletov2010} and to incorporate
the noise intensity, $D$, into dynamics in accordance with \prettyref{eq:Dyn-D-rate}.
This procedure is just reasonable since involves dynamical fluctuations
of the membrane electrical characteristics. To correctly sample the
invariant measure, dynamics must be ergodic.

Consider dynamics in the extended phase space $\left(\left\{ u_{i}\right\} ,\lambda,\left\{ \eta_{i}\right\} \right)$,
\begin{equation}
\begin{array}{l}
\dot{u_{i}}=\end{array}-\lambda\partial_{i}\mathcal{V}[u]+\eta_{i},\quad\dot{\lambda}=g(u),\quad\dot{\eta_{i}}=h_{i}(u),\quad i\in\mathbb{Z};\label{eq:ND}\end{equation}
functions $g(u)$ and $h_{i}(u)$ are to be determined. The extra
dynamical variables $\lambda$ and $\eta_{i}$ model the environment
and thus they represent the noise effect on the Nagumo dynamics.
\begin{rem*}
Term $\eta_{i}$ in the dynamical equations (\ref{eq:ND}) is important.
Indeed, assume $\eta_{i}\equiv0$. Then, at an equilibrium $\partial_{i}\mathcal{V}[u]=0$,
the evolution comes to halt and no longer fluctuates, irrespective
of the time dependence of $\lambda$. For initial conditions with
$\partial_{i}\mathcal{V}[u]\neq0$ after a time variable rescaling,
it is a gradient flow as defined in \citep{Katok1996}, and all phase
space trajectories moves along paths with equilibrium points at either
end. Thus dynamics is not ergodic. For a further discussion we refer
to \citep{SDC07,Samoletov2010}.
\end{rem*}
To determine functions $g(u)$ and $h_{i}(u)$, calculate, on the
analogy of \prettyref{eq:SLyap},\begin{equation}
\dot{\mathcal{V}}[u]=-\lambda\sum_{\{i\}}\left(\partial_{i}\mathcal{V}[u]\right)^{2}+\sum_{(i)}\eta_{i}\partial_{i}\mathcal{V}[u].\label{eq:DynFluct}\end{equation}
Respect to the second term on r.h.s. of \prettyref{eq:DynFluct} we
put the following requirement to the time average, $\overline{\sum_{(i)}\eta_{i}\partial_{i}\mathcal{V}[u]}=0$.
A series of $\eta$-dynamics satisfies this condition. Two principal
limit cases are: fluctuations of current in different nodes are independent
or synchronous. Correspondingly we endow variables $\left\{ \eta_{i}\right\} $
with the following dynamical equations,\begin{equation}
\dot{\eta}_{i}\sim\partial_{i}\mathcal{V}[u],\quad i\in\mathbb{Z},\quad\textrm{and}\quad\dot{\eta}_{i}\sim\sum_{(j)}\partial_{j}\mathcal{V}[u],\quad\forall i\in\mathbb{Z}.\label{eq:eta}\end{equation}
However, respect to the first term in r.h.s. of \prettyref{eq:DynFluct},
we cannot repeat the trick and set $\dot{\lambda}\sim\sum_{\{i\}}\left(\partial_{i}\mathcal{V}[u]\right)^{2}$,
since this results in no noise effect. To overcome this difficulty,
we implement conjecture \prettyref{eq:Dyn-D-rate} into $\lambda$-dynamics
and \prettyref{eq:DynFluct}, and explicitly set\begin{equation}
\dot{\lambda}\sim\sum_{\{i\}}\left(\partial_{i}\mathcal{V}[u]\right)^{2}-D\sum_{\{i\}}\partial_{i}^{2}\mathcal{V}[u].\label{eq:lambda-da}\end{equation}

\begin{lem}
Assume $\lambda$ to be bounded variable, its dynamics is given by
\prettyref{eq:lambda-da} and $\overline{\sum_{(i)}\eta_{i}\partial_{i}\mathcal{V}[u]}=0$
(e.g. one of dynamical equations \prettyref{eq:eta}). Then $\overline{\lambda\sum_{\{i\}}\partial_{i}^{2}\mathcal{V}[u]}=0.$ \end{lem}
\begin{proof}
First we multiply \prettyref{eq:lambda-da} by $\lambda$ and take
into account equation \prettyref{eq:DynFluct}. Then we apply the
time averaging to the resulted equation. Thus we easily accomplish
lemma. Indeed,\[
0=\overline{-\lambda\left[\sum_{\{i\}}\left(\partial_{i}\mathcal{V}[u]\right)^{2}-D\sum_{\{i\}}\partial_{i}^{2}\mathcal{V}[u]\right]}=\overline{-\dot{\mathcal{V}}[u]+\sum_{(i)}\eta_{i}\partial_{i}\mathcal{V}[u]-D\lambda\sum_{\{i\}}\partial_{i}^{2}\mathcal{V}[u]}=D\overline{\lambda\sum_{\{i\}}\partial_{i}^{2}\mathcal{V}[u]}.\]

\end{proof}
This lemma together with the equations \prettyref{eq:DynFluct}-\prettyref{eq:lambda-da}
allows us to determine functions $g(u)$ and $h_{i}(u)$ explicitly,\begin{equation}
g=\frac{1}{Q_{\lambda}}\sum_{(i)}\left[(\partial_{i}\mathcal{V}[u])^{2}-D\partial_{i}^{2}\mathcal{V}[u]\right],\quad h_{i}=-\frac{1}{Q_{\eta}}\partial_{i}\mathcal{V}[u]\;\mathrm{or}\; h_{i}=-\frac{1}{Q_{\eta}}{\displaystyle \sum_{(j)}}\partial_{j}\mathcal{V}[u],\quad i\in\mathbb{Z},\label{eq:g-h}\end{equation}
where $Q_{\lambda}$ and $Q_{\eta}$ are parameters. Variables $\eta_{i}$
and corresponding functions $h_{i}$ and not unique and dynamical
equations can be simplified.

To verify the requirement on the invariant measure we prove the theorem.
\begin{thm}
Assume the extended dynamics in form \prettyref{eq:ND} where functions
$g$ and $h_{i}$ are given by \prettyref{eq:g-h}, $Q_{\lambda}>0$
and $Q_{\eta}>0$. Then the augmented measure, \begin{equation}
d\mu\sim exp\{-D^{-1}\mathcal{V}[u]\}exp[\{-D^{-1}(\frac{1}{2}Q_{\lambda}\lambda^{2}+\frac{1}{2}Q_{\eta}\sum_{(i)}\eta_{i}^{2})\}\prod_{(i)}du_{i}d\lambda d\eta_{i}=\rho_{\infty}\prod_{(i)}du_{i}d\lambda d\eta_{i},\label{eq:IMV}\end{equation}
is invariant for the extended dynamics\emph{.}\end{thm}
\begin{rem*}
It should be noted that the $\eta$-dynamics is not unique and so
correspondingly it allows a variety of $\eta$-factors of the augmented
measure, although they are still Gaussian. \emph{E.g}., with the synchronous
dynamical fluctuations, $\dot{\eta}=-\frac{1}{Q_{\eta}}\sum_{(j)}\partial_{j}\mathcal{V}[u]$,
we arrive at $exp[\{-D^{-1}\frac{1}{2}Q_{\eta}\eta^{2}\}d\eta$. However
all cases can be treated analogously.\end{rem*}
\begin{proof}
The Liouville operator corresponding to the dynamics in the extended
phase space (\ref{eq:ND}) has the form, $\mathcal{L}^{*}\rho=-\sum_{(i)}\partial_{i}\left[(-\lambda\partial_{i}\mathcal{V}[u]+\eta_{i})\rho\right]-\partial_{\lambda}[g(u)\rho]-\sum_{(i)}\partial_{\eta_{i}}[h_{i}(u)\rho],$
and the Liouville equation reads $\partial_{t}\rho=\mathcal{L}^{*}\rho$.
Therefore, to prove the theorem we have to prove the identity, $\mathcal{L}^{*}\rho_{\infty}\equiv0$,
that means the dynamics \prettyref{eq:ND} preserves the augmented
measure (\ref{eq:IMV}). A straightforward calculation of all partial
derivatives that are involved in $\mathcal{L}^{*}\rho_{\infty}$ with
further simplification brings to the required identity, $\mathcal{L}^{*}\rho_{\infty}\equiv0$.
The theorem is proved. 
\end{proof}
From the perspective of numerical simulations and further mathematical
analysis, \emph{e.g.} the Hamiltonian representation of the proposed
dynamics, it is important to find a first integral of motion. We accomplish
this task with the following lemma.
\begin{lem}
Let the dynamical system \prettyref{eq:ND} and \prettyref{eq:g-h}
be augmented with the redundant dynamical variable $\zeta$, $\dot{\zeta}=-\lambda\sum_{\{i\}}\partial_{i}^{2}\mathcal{V}[u]$.
Then $I=\mathcal{V}[u]+\frac{1}{2}Q_{\lambda}\lambda^{2}+\frac{1}{2}Q_{\eta}\sum_{(i)}\eta_{i}^{2}-D\zeta$
is the first integral of the augmented dynamical system.\end{lem}
\begin{proof}
We derive $\dot{I}=0$ by direct calculation.\end{proof}
\begin{rem*}
Since the origin of coordinates of the redundant variable $\zeta$
is arbitrary, it is always possible for an arbitrary fixed trajectory
to set $I=0$. $I$ is apparent control parameter in numerical simulations.
Besides, $I$ is related to $\rho_{\infty}$ and thus can be considered
from a perspective of the Hamiltonian reformulation of dynamics on
the level set $I=0$ \citep{SDC07}. However we do not consider this
problem here.
\end{rem*}
We can now ask whether the dynamics \prettyref{eq:ND},\prettyref{eq:g-h}
is ergodic. There is no a definite answer to this question. Following
\citep{SDC07} we can apply the Frobenius theorem of differential
geometry \citep{lang2002introduction} but this provides with a partial
answer only. Here, in order to provide ergodicity, we adopt the method
proposed in \citep{SDC07} and rigorously investigated in \citep{Leimkuhler2009}.
Namely, we add a Gaussian random noise to the $\lambda$-dynamics.
In contrast to the model outlined in \prettyref{sec:Random-noise},
where stochastic currents are added at each node, this approach relies
on single and indirect stochastic perturbation. Experiments \citep{Leimkuhler2009}
reveal that, in context of the molecular dynamics, it results in a
relatively weak perturbance effect on deterministic dynamics. Thus,
we reformulate $\lambda$-dynamics \prettyref{eq:ND} in the form,\begin{equation}
\begin{array}{l}
\dot{u_{i}}=\end{array}-\lambda\partial_{i}\mathcal{V}[u]+\eta_{i},\quad\dot{\lambda}=g(u)-\gamma\lambda+\sqrt{2\gamma DQ_{\lambda}^{-1}}\xi(t),\quad\dot{\eta_{i}}=h_{i}(u),\quad i\in\mathbb{Z},\label{eq:NDstoch}\end{equation}
where $\gamma>0$ is a parameter.
\begin{thm}
Assume stochastically perturbed extended dynamics in the form \prettyref{eq:NDstoch}
where functions $g$ and $h$ are given by \prettyref{eq:g-h}, $Q_{\lambda}>0$,
$Q_{\eta}>0$. Then the augmented measure \prettyref{eq:IMV} is invariant
for this dynamics\emph{.}\end{thm}
\begin{proof}
The Fokker-Planck operator corresponding to (\ref{eq:NDstoch}) has
the form, $\mathcal{F}^{*}\rho=\mathcal{L}^{*}\rho+\gamma\partial_{\lambda}\left[\left(\lambda+DQ_{\lambda}^{-1}\partial_{\lambda}\right)\rho\right],$
and the Fokker-Planck equation reads $\partial_{t}\rho=\mathcal{F}^{*}\rho$.
After a series of routine calculations we arrive at $\mathcal{F}^{*}\rho_{\infty}\equiv0$.
Thus the stochastically perturbed dynamics \prettyref{eq:NDstoch}
preserves the augmented measure (\ref{eq:IMV}). 
\end{proof}

\subsection*{Test simulations. Single cell dynamics. }

Low dimensional systems often reveal the ergodicity problem in a probability
distribution dynamical sampling. For this reason, it is important to
test the presented noise modelling method capable of generating the
right statistic for a single Nagumo cell. We choose for this purpose
$f(u)=-4u(u-\alpha)(u-1)$. Simulations are performed using global
parameters $D=0.04$ and $\gamma=1$, for $t=10^{6}$.
\begin{figure}[h]
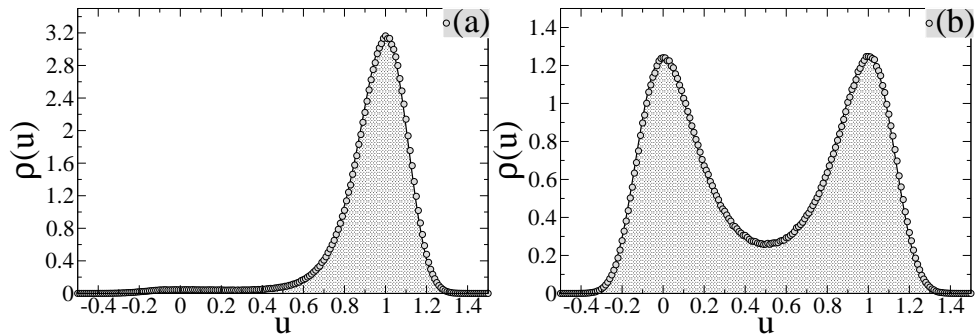

\centering{}\includegraphics[clip,scale=0.25]{a.eps}\includegraphics[clip,scale=0.25]{b.eps}\caption{Probability distributions of $u$ variable (shown on background of
exact analytical distribution). Densities are calculated as normalized
sojourn distributions. Correspondingly, $(\mathrm{a})$ $\alpha=0.25$
and $(\mathrm{b})$ $\alpha=0.5$.}
\end{figure} 
Figure 1 shows the probability distribution of the $u$ variable calculated with
the dynamical equations and compared with exact analytical distribution.
Their solid agreement brings a severe test of our approach.

\section{Conclusion}

We have presented a novel mathematical approach to model noise in
dynamical systems. We do so by considering dynamics of a chain of
diffusively coupled Nagumo cells affected by noise. We have shown
that the noise in transmembrane current can be effectively modelled
as fluctuations in electric characteristics of the membrane. Test
simulations give a solid support to the mathematical scheme. The proposed
approach to model noise in a nerve fibre is different from the standard
additive stochastic current perturbation and thus demonstrates a potential
for further application.

\section*{Acknowledgments }

This work was supported in part by the University of Liverpool. AS would like to thank 
Department of Mathematical Sciences for hospitality.

\end{document}